\documentclass{amsart}
\numberwithin{equation}{section}
\usepackage{amscd}
\usepackage[all]{xy}
\CompileMatrices
\usepackage{amssymb}
\let\cal\mathcal
\def\Ascr{{\cal A}}

\def\Cscr{{\cal C}}

\def\Fscr{{\cal F}}

\def\Tscr{{\cal T}}

\let\blb\mathbb

\def \ZZ{{\blb Z}}

\def \NN{{\blb N}}

\def\Id{\operatorname{id}}

\def\Mod{\operatorname{Mod}}
\def\mod{\operatorname{mod}}
\def\Gr{\operatorname{Gr}}

\def\gr{\operatorname{gr}}

\def\gr{\operatorname {gr}}

\def\Ext{\operatorname {Ext}}
\def\Hom{\operatorname {Hom}}

\def\im{\operatorname {im}}

\def\Tr{\operatorname {Tr}}

\def\rk{\operatorname {rk}}

\def\r{\rightarrow}

\DeclareMathOperator{\Tors}{Tors}

\let\dirlim\injlim

\newtheorem{lemma}{Lemma}[section]
\newtheorem{proposition}[lemma]{Proposition}
\newtheorem{theorem}[lemma]{Theorem}

\theoremstyle{definition}

{

}

\theoremstyle{remark}

\newtheorem{remark}[lemma]{Remark}

\newdimen\uboxsep \uboxsep=1ex
\def\uboxn#1{\vtop to 0pt{\hrule height 0pt depth 0pt\vskip\uboxsep
\hbox to 0pt{\hss #1\hss}\vss}}

\def\uboxs#1{\vbox to 0pt{\vss\hbox to 0pt{\hss #1\hss}
\vskip\uboxsep\hrule height 0pt depth 0pt}}

\title{On the {\boldmath $\ZZ D_\infty$} category}
\author{Michel Van den Bergh}
 \address{Departement WNI,  Limburgs Universitair Centrum, 
 3590 Diepenbeek, Belgium.}
  \email{vdbergh@luc.ac.be}
\thanks{The author is a senior researcher at the FWO}

\begin{document}
\begin{abstract}
In this paper we give a direct proof of the
properties of the $\ZZ D_\infty$ category which was
introduced in the classification of noetherian,
hereditary categories with Serre duality by Idun Reiten
and the author.
\end{abstract}
\maketitle
\section{Introduction}
Below $k$ is a field. All categories will be $k$-linear. An abelian or
triangulated category $\Ascr$ is  $\Ext$-finite if for all
objects $A,B\in \Ascr$ one has that $\oplus_i \Ext^i(A,B)$ is finite
dimensional. If $\Ascr$ is triangulated and $\Ext$-finite then we say
that $\Ascr$ satisfies Serre duality \cite{Bondal4} if there exists an
auto-equivalence $F$ of $\Ascr$ together with  isomorphisms
\[
\Hom_\Ascr(A,B)\r \Hom_\Ascr(B,FA)^\ast
\]
natural in $A$, $B$ 
(where $(-)^\ast$ is the $k$-dual).
If $\Ascr$ is abelian and $\Ext$-finite then we say that $\Ascr$ satisfies
Serre duality if this is the case for $D^b(\Ascr)$. The following
result can be extracted from \cite[Ch.\ 1]{ReVdB1}.
\begin{theorem} 
\label{ref-1.1-0}
Assume that $\Cscr$ is an $\Ext$-finite hereditary
category without injectives or projectives. Then the following are 
equivalent 
\begin{enumerate}
\item $\Cscr$ has almost split sequences.
\item $\Cscr$ satisfies Serre duality.
\item There is an auto-equivalence $V:\Cscr\r \Cscr$ together with
natural isomorphisms
\begin{equation}
\label{ref-1.1-1}
\Hom_\Cscr(A,B)\r \Ext^1_\Cscr(B,VA)^\ast
\end{equation}
\end{enumerate}
Furthermore the functor $V$  coincides with the Auslander-Reiten
translate $\tau$ when evaluated on objects. 
\end{theorem}

In the classification of noetherian $\Ext$-finite hereditary
categories with Serre duality in \cite{ReVdB1} we considered a category
$\Cscr$ defined by the following pullback diagram
\begin{equation}
\label{ref-1.2-2}
\begin{CD}
\mod(k)\oplus \mod(k) @>>> \mod(k)\\ @AAA @AAA\\\Cscr @>>> \gr(k[x])
\end{CD}
\end{equation}
where the horizontal map sends $(V_1,V_2)$ to $V_1\oplus V_2$ and the
vertical map is localizing at $x$ followed by restricting to degree zero. It
was shown by a rather indirect argument that $\Cscr$ is a noetherian,
$\Ext$-finite, hereditary abelian category without injectives or
projectives which satisfies Serre duality.  It was also shown
that the AR-quiver of $\Cscr$ has two components, one equal to $\ZZ A_\infty$
(a ``wing'') and the other equal to $\ZZ D_\infty$. For this reason
$\Cscr$ was called the ``$\ZZ D_\infty$-category''.

The aim of this paper is to give a direct proof of the above facts. In
addition we will also establish a link with one-dimensional $A_n$
singularities.
\section{Elementary properties}
It is easy to see that $\Cscr$, as defined in the introduction is a
noetherian abelian category.  It will be convenient to consider the
locally noetherian Grothendieck category $\tilde{\Cscr}$ associated to
$\Cscr$. It follows for example by \cite[Prop.\ 2.14]{lowenvdb1} that $D^b(\Cscr)$ and
$D^b_\Cscr(\tilde{\Cscr})$ are equivalent.  Hence the
$\Ext$-groups between objects in $\Cscr$ may be computed in $\tilde{\Cscr}$.

The objects of $\tilde{\Cscr}$
are quadruples $(M,V_0,V_1,\phi)$ where $M$ is a graded $k[x]$-module,
 $V_0$, $V_1$ are $k$-vector spaces and $\phi$ is an isomorphism of
 $k[x]$ modules $M_x\r (V_0\oplus V_1)\otimes_k k[x,x^{-1}]$. Objects in $\Cscr$ are given by the quadruples
$(M,V_0,V_1,\phi)$ in which $M$ is finitely generated.

Sending $(M,V_0,V_1,\phi)$ to $M$ defines an faithful exact functor
$\tilde{\Cscr}\r \Gr(k[x])$ which we call the restriction functor and
which we denote by $(-)_{k[x]}$. 

We write
$M(n)=(M(n),V_0,V_1,\phi)$ where we have identified $M(n)_x$ with
$M_x$ through multiplication with $x^n$. Furthermore we define $\sigma
(M)=(M,V_1,V_0,\phi)$. 

We define
$\tilde{\Tscr}\subset\tilde{\Cscr}$ and $\tilde{\Fscr}\subset
\tilde{\Cscr}$ respectively as the inverse images of the $x$-torsion
and $x$-torsion free modules in $\Gr(k[x])$.  $\Tscr$ and $\Fscr$ are
defined similarly, but starting from $\Cscr$.

By $\tilde{\Cscr}_x$ we denote the full subcategory of $\tilde{\Cscr}$
with objects the quadruples $(M,V_0,V_1,\phi)$ in which $x$ acts
invertibly on $M$. 

 We denote by $(-)_x$ the functor $\tilde{\Cscr}\r
\tilde{\Cscr}_x$ which sends $(M,V_0,V_1,\phi)$ to
$(M_x,V_0,V_1,\phi)$. Clearly if $M\in\tilde{\Cscr}$ and
$N\in\tilde{\Cscr}_x$ then the canonical maps
\begin{equation}
\label{ref-2.1-3}
\Hom_{\tilde{\Cscr}}(M,N)\r \Hom_{\tilde{\Cscr}_x}(M,N)\r\Hom_{{\tilde{\Cscr}_x}}(M_x,N)
\end{equation}
are isomorphisms. We list a few other obvious facts.
\begin{itemize}
\item[(O1)] 
\label{ref-2-4}
$(\tilde{\Tscr},\tilde{\Fscr})$ forms a torsion pair in
  $\tilde{\Cscr}$. That is $\Hom(\tilde{\Tscr},\tilde{\Fscr})=0$ and
  for any $M\in\tilde{\Cscr}$ there exists an exact sequence
  (necessarily unique)
\[
0\r T\r M\r F\r 0
\]
with $T\in \tilde{\Tscr}$ and $F\in \tilde{\Fscr}$. 
\item[(O2)] 
\label{ref-2-5}
If $T\in \tilde{\Tscr}$ and $M\in\tilde{\Cscr}$ then 
\begin{align*}
\Hom_{\tilde{\Cscr}}(T,M)&=\Hom_{k[x]}(T,M)\\
\Hom_{\tilde{\Cscr}}(M,T)&=\Hom_{k[x]}(M,T)
\end{align*}
\item[(O3)]
\label{ref-2-6}
The restriction functor defines an equivalence between $\tilde{\Tscr}$
and $\Tors(k[x])$ where $\Tors(k[x])$ denotes the $x$-torsion modules
is $\Gr(k[x])$. 
\item[(O4)] 
\label{ref-2-7}
The functor $\tilde{\Cscr}_x\r \Mod(k)\oplus \Mod(k)$ which sends
$(M,V_0,V_1,\phi)$ to $V_0\oplus V_1$ is an equivalence of categories.
\end{itemize}
\def\Oref#1{(O\ref{#1})}
Combining (O4) with \eqref{ref-2.1-3} yields in particular 
\begin{itemize}
\setcounter{enumi}{4}
\item[(O5)]
\label{ref-2-8}
If $N\in \tilde{\Cscr}_x$ then $\Hom_{\Cscr}(-,N)$ is exact. Hence the
objects in $\tilde{\Cscr}_x$ are injective in $\tilde{\Cscr}$.
\end{itemize}
We now describe the indecomposable injectives in $\tilde{\Cscr}$. For
$n\in \ZZ$ let $E_{n}$ be the graded injective $k[x]$-module
given by $k[x,x^{-1}]/x^{n+1} k[x]$. Since $E_n\in \Tors(k[x])$ there
exists by (O3) a corresponding object in $\tilde{\Tscr}$
which we denote by the same symbol.
From (O2) it follows
that $\Hom(-,E_{n})$ is exact and hence $E_{n}$ is injective
in $\tilde{\Cscr}$.

To construct other injectives we note that by (O5) we know that the
objects in $\tilde{\Cscr}_x$ are injective in $\tilde{\Cscr}$. Since
by (O4) $\tilde{\Cscr}_x$ is equivalent to $\Mod(k)\oplus \Mod(k)$ there must
be two corresponding indecomposable injectives in $\tilde{\Cscr}$.
They are given by
\begin{align*}
E^0&=(k[x,x^{-1}],k,0,\Id_{k[x,x^{-1}]})\\
E^1&=(k[x,x^{-1}],0,k,\Id_{k[x,x^{-1}]})
\end{align*}
\begin{proposition}
\label{ref-2.1-9}
\begin{enumerate}
\item
 $(E_{ n})_n,E^0,E^1$ forms a complete list of
indecomposable injectives in $\tilde{\Cscr}$.
\item Every object in $\tilde{\Cscr}$ has injective dimension one (and hence
$\tilde{\Cscr}$ and $\Cscr$ are hereditary \cite[Prop.\ A.3]{ReVdB1}).
\item $\Cscr$ is $\Ext$-finite.
\end{enumerate}
\end{proposition}
\begin{proof}
Since the listed injectives are clearly indecomposable (1) follows
if we can show that any indecomposable object can be embedded in a direct sum
of them \cite{Matlis}. 

  We prove (1) and (2) together by showing that every object
  $M\in\tilde{\Cscr}$ has a resolution of length at most two whose terms
  consist of direct sums of the injectives given in (1). 
 By
  (O1) it is clearly sufficient to prove this claim
  separately in the cases $M\in\tilde{\Tscr}$ and $M\in\tilde{\Fscr}$.

Assume first that $M\in \tilde{\Tscr}$. Then $M$ has an injective
resolution
\begin{equation}
\label{ref-2.2-10}
0\r M \r I_0\r I_1\r 0
\end{equation}
in $\Tors(k[x])$. By (O3) this resolution corresponds to one in
$\tilde{\Cscr}$. Furthermore by the structure of the injectives in
$\Gr(k[x])$ the $I_i$ are direct sums of the $E_{ n}$ in
$\Gr(k[x])$. Again by (O3) the same is true in
$\tilde{\Cscr}$.

Now assume that $M\in\tilde{\Fscr}$. Consider the short exact
sequence
\begin{equation}
\label{ref-2.3-11}
0\r M\r M_x\r M_x/M\r 0
\end{equation}
$M_x$ lies in $\tilde{\Cscr}_x$ and hence by (O4) is a direct sum of
copies of $E^0$ and $E^1$. $M_x/M$ is $x$-divisible and lies in
$\tilde{\Tscr}$ and so by (O3) $M_x/M$ is a direct sum of copies of
$E_{ n}$. Whence \eqref{ref-2.3-11} is the kind of resolution
we were looking for.

To prove (3) we note that if $E,F$ are indecomposable injectives as
in (1) then $\dim \Hom_{\tilde{\Cscr}}(E,F)\le 1$. Thus it suffices to show
that every $M\in\Cscr$ has an injective resolution consisting in every degree of
a finite number of indecomposable injectives.  This follows
easily from the construction.
\end{proof}
\begin{proposition} 
\label{ref-2.2-12} If $F\in \tilde{\Fscr}$ and $T\in \tilde{\Tscr}$ then
$\Ext^1_{\tilde{\Cscr}}(F,T)=0$. In particular every object in $\Cscr$ is 
of the form $F\oplus T$ with $F\in \Fscr$ and $T\in \Tscr$.
\end{proposition}
\begin{proof}
  It follows from (O2) that $\Hom(F,-)$ is exact on $\Tscr$. Since by
  the proof of the previous proposition $T$ has a $\tilde{\Cscr}$
  injective resolution inside $\tilde{\Tscr}$, we are done.
\end{proof}
Now we describe the $\Ext$-groups between objects in
$\tilde{\Fscr}$.
\begin{lemma} 
\label{ref-2.3-13}
Assume that $F=(F,V_0,V_1,\phi)$,
  $F'=(F',V'_0,V'_1,\phi')$ are objects in $\tilde{\Fscr}$. Then there exists
  an exact sequence of the form 
\begin{equation}
\label{ref-2.4-14}
0\r \Hom_{\tilde{\Cscr}}(F,F')\r \Hom_{k[x]}(F,F')\r 
\Hom_k(V_0,V'_1)\oplus \Hom_k(V_1,V'_0)\r \Ext^1_{\tilde{\Cscr}}(F,F')\r 0
\end{equation}
\end{lemma}
\begin{proof}
We start with the short exact sequence
\[
0\r F'\r F'_x\r F'_x/F'\r 0
\]
which according to the proof of lemma \ref{ref-2.1-9}
is an injective resolution of $F'$, both in $\tilde{\Cscr}$ and in
$\Gr(k[x])$. 

Applying $\Hom_{\Cscr}(F,-)$, $\Hom_{\Gr(k[x])}(F,-)$ and comparing
yields a commutative diagram with exact rows and columns.
{\tiny
$$\hskip -1cm\xymatrix{
& & 0\ar[d]\\
0
\ar[r]& \Hom_{\tilde{\Cscr}}(F,F')\ar[r]\ar[d]& \Hom_k(V_0,V'_0)\oplus
\Hom_k(V_1,V_1') \ar[r]\ar[d]& 
\Hom_{\tilde{\Cscr}}(F,F'_x/F')\ar[r]\ar@{=}[d]& \Ext^1_{\tilde{\Cscr}}(F,F') \ar[r] &0\\
0 \ar[r] &\Hom_{\Gr(k[x])}(F,F')\ar[r]& \Hom_k(V_0\oplus V_1,V'_0\oplus V'_1)
\ar[r]\ar[d]&
\Hom_{\Gr(k[x])}(F,F'_x/F') \ar[r]& 0\\
&& \Hom(V_0,V'_1)\oplus \Hom(V_1,V'_0)\ar[d]\\
&& 0
}$$
}
\eqref{ref-2.4-14} now follows from the previous diagram
through an easy diagram chase.
\end{proof}
\begin{proposition} $\Cscr$ has neither injectives nor projectives. 
\end{proposition}
\begin{proof}
Since $\Tscr$ is equivalent to the $x$-torsion modules
in $\gr(k[x])$,
it is easy to see that $\Tscr$ does not contain any injective or
projectives.

If $0\neq (F,V_0,V_1,\phi)$ in $\Fscr$ then by considering the
faithful restriction functor to $\gr(k[x])$ we see that
$\Hom_\Cscr(F,\sigma F(-n))=0$ for $n\gg 0$. On the other hand $V_0$ or
$V_1\neq 0$. It follows from the previous lemma that $\Ext^1_\Cscr(F,\sigma F(-n))\neq 0$ for $n\gg 0$. Hence $F$ is not projective. A similar
argument shows that $F$ is not injective.
\end{proof}
\begin{remark}
The reason why we called this section ``Elementary properties'' is that
the stated results   hold in greater generality. For example,
suitably adapted versions would be  valid for the pullback of
\[
\begin{CD}
\mod(k)^{\oplus m} @>>> \mod(k)\\ @. @AAA\\ @. \gr(k[x])
\end{CD}
\]
for any $m$. By contrast, the results in the next section require
$m=2$.
\end{remark}
\section{Serre duality}
\S\label{ref-3-15}
Our next aim is to prove that $\Cscr$ satisfies Serre duality. First
we construct a Serre functor on $\Fscr$. Put $VM=\sigma(M)(-1)$.

The first step in proving Serre duality is constructing a ``trace
map'' $\eta_M:\Ext^1_\Cscr(M,VM)\r k$ for $M\in\Fscr$ which should
corresponds to the identity map in $\Hom_\Cscr(M,M)$ under the isomorphism \eqref{ref-1.1-1}.

We now use \eqref{ref-2.4-14} to construct the trace map $\eta_F$
for $F=(F,V_0,V_1,\phi)\in\Fscr$. In this case
$VF=(F(-1),V_1,V_0,\phi)$ and we have an exact sequence
\[
\Hom_{k[x]}(F,VF)\r \Hom_k(V_0,V_0)\oplus \Hom_k(V_1,V_1)\r
\Ext^1_{\Cscr}(F,VF)\r  0
\]
\begin{lemma}
\label{ref-3.1-16}
The composition 
\begin{equation}
\label{ref-3.1-17}
\Hom_{k[x]}(F,VF)\r \Hom_k(V_0,V_0)\oplus
\Hom_k(V_1,V_1)\xrightarrow{\Tr_{V_0}+\Tr_{V_1}} k
\end{equation}
is the zero map. 
\end{lemma}
\begin{proof}
To see this note that 
$\Hom_{k[x]}(F,VF)=\Hom(F,F(-1))$ and furthermore that
\eqref{ref-3.1-17} can be extended to a commutative diagram.
$$\xymatrix{
\Hom_{k[x]}(F,F(-1))\ar[r]& \Hom_k(V_0\oplus V_1,V_0\oplus
V_1)\ar[r]\ar[d]^{\Tr} 
 & \Hom_k(V_0,V_0)\oplus
\Hom_k(V_1,V_1)\ar[ld]^{\Tr_{V_0}+\Tr_{V_1}}\\
&k
}
$$
By choosing a basis for $F$ as graded $k[x]$-module one easily sees
that every element of $\Hom(F,F(-1))\subset \Hom(F,F)$ is
nilpotent. Since nilpotent elements have zero trace it follows that
the composition
\[
\Hom_{k[x]}(F,F(-1))\r \Hom_k(V_0\oplus V_1,V_0\oplus
V_1)\xrightarrow{\Tr} k
\]
is zero. This proves what we want.
\end{proof}
From lemma \ref{ref-3.1-16} together with
\eqref{ref-2.4-14} there exists a unique map
$\eta_F:\Ext^{1}_{\Cscr}(F,VF)\r k$ which makes the following diagram
commutative.
$$\xymatrix{
\Hom_{k[x]}(F,VF)\ar[r]& \Hom(V_0,V_0)\oplus \Hom(V_1,V_1)
\ar[r]\ar[dr]^{\Tr_{V_0}+\Tr_{V_1}}& \Ext^1_{\Cscr}(F,VF)\ar[r]
\ar[d]^{\eta_F} & 0\\
&&k
}
$$
To continue it will be convenient to use the Yoneda multiplication 
on $\Ext^\ast_\Cscr(-,-)$. In order
to have compatibility with the notation for compositions of maps we
will 
write the Yoneda multiplication as a pairing
\[
\Ext^\ast_\Cscr(B,C)\times \Ext^\ast_\Cscr(A,B)\r \Ext^\ast_\Cscr(A,C)
\]
We extend $\eta_F$ to a map $\Ext^\ast(F,VF)\r k$ by letting it act
trivially on $\Hom(F,VF)$. 
\begin{lemma}
\label{ref-3.2-18}
 Let $F,G\in \Fscr$ and assume that $f\in
  \Ext^\ast_{\Cscr}(F,G)$ and $g\in \Ext^\ast_{\Cscr}(G,VF)$. Then we
  have $\eta_F (gf)=\eta_G(V(f)g)$.
\end{lemma}
\begin{proof}
We may assume that $f$ and $g$ are homogeneous. Furthermore the cases
where $f$,$g$ are both of degree $0$ or of degree $1$ are
trivial. Hence we may assume that $(\deg f,\deg g)=(0,1)$ or $(\deg
f,\deg g)=(1,0)$.

Let us consider the first possibility. We check that
$\eta_F(-f)=\eta_G(V(f)-)$ as maps  $\Ext^1(G,VF)\r k$.
This amounts to the commutativity of
\begin{equation}
\label{ref-3.2-19}
\begin{CD}
\Ext^1(G,VF) @>>> \Ext^1(G,VG)\\
@VVV @V\eta_G VV\\
\Ext^1(F,VF) @>\eta_F>> k
\end{CD}
\end{equation}
Assume that $F=(F,V_0,V_1,\phi)$, $G=(G,W_0,W_1,\theta)$. Then $f$ induces
maps $f_0:V_0\r W_0$ and $f_1:V_1\r W_1$.
Elementary linear algebra yields that we
  have a commutative diagram
\[
\begin{CD}
\Hom(W_0,V_0)\oplus \Hom(W_1,V_1)@>>> \Hom(W_0,W_0)\oplus
\Hom(W_1,W_1) \\
@VVV @V\Tr_{W_0}+\Tr_{W_1} VV\\
\Hom(V_0,V_0)\oplus
\Hom(V_1,V_1) @>\Tr_{V_0}+\Tr_{V_1}>> k
\end{CD}
\]
This diagram, together with the definition of $\eta$ yields the
commutativity of \eqref{ref-3.2-19}.

Now we consider the possibility $(\deg
f,\deg g)=(1,0)$. Since we trivially have 
\begin{equation}
\label{ref-3.3-20}
\eta_{VX}\circ V=\eta_X  
\end{equation}
it is sufficient to prove that $\eta_{VF}
(V(g)V(f))=\eta_G(V(f)g)$. Replacing $(Vf,g)$ by $(g,f)$ this reduces
to the previous case.
\end{proof}
We are now in a position to prove Serre duality for objects in
$\Fscr$. We will show that the pairing
\begin{equation}
\Hom_{\Cscr}(F,G)\times \Ext^1_{\Cscr}(G,VF)\r
\Ext^1(F,VF)\xrightarrow{\eta_F} k:(f,g)\mapsto \eta_F (gf)
 \label{ref-3.4-21}
\end{equation}
is non-degenerate. By lemma \ref{ref-3.2-18} the non-degeneracy
of \eqref{ref-3.4-21} for all $F,G$ is equivalent to the non-degeneracy
of the pairing
\begin{equation}
\Ext^1_{\Cscr}(F,G)\times \Hom_{\Cscr}(G,VF)\r
\Ext^1(F,VF)\xrightarrow{\eta_F} k:(f,g)\r \eta_F(gf)
 \label{ref-3.5-22}
\end{equation}
for all $F,G$. It follows also easily from lemma
\ref{ref-3.2-18} that \eqref{ref-3.4-21} and \eqref{ref-3.5-22}
are natural in $F$ and $G$.
\begin{lemma}
\label{ref-3.3-23}
If we have and exact sequence
\begin{equation}
\label{ref-3.6-24}
0\r F_1\r F\r F_2\r 0
\end{equation}
in $\Fscr$ and if we have non-degeneracy of \eqref{ref-3.4-21} and
\eqref{ref-3.5-22} for two out of the three pairs $(F_1,G), (F,G),
(F_2,G)$ then we also have it for the third one.  A similar statement
holds for an exact sequence
\begin{equation}
\label{ref-3.7-25}
0\r G_1\r G\r G_2\r 0
\end{equation}
\end{lemma}
\begin{proof} Assume that we have an exact sequence of the form
  \eqref{ref-3.6-24}.  We claim that the following diagram
  with exact rows is commutative. 
 {\tiny
$$\hskip -1cm\xymatrix@-0.5cm{
0 \ar[r]&
\Hom(F_2,G) \ar[r] \ar[d]^{\alpha_2}&
\Hom(F,G) \ar[r]\ar[d]^{\alpha} &
\Hom(F_1,G) \ar[r]\ar[d]^{\alpha_1} &
\Ext^1(F_2,G) \ar[r]\ar[d]^{\beta_2} &
\Ext^1(F,G) \ar[r] \ar[d]^{\beta}&
\Ext^1(F_1,G) \ar[r] \ar[d]^{\beta_1}&
0\\
0 \ar[r]&
\Ext^1(G,VF_2)^\ast \ar[r]&
\Ext^1(G,VF)^\ast \ar[r] &
\Ext^1(G,VF_1)^\ast \ar[r] &
\Hom(G,VF_2)^\ast \ar[r] &
\Hom(G,VF)^\ast \ar[r]&
\Hom(G,VF_1)^\ast \ar[r]&
0
}
$$
}
Here the maps labeled by $\alpha$ are obtained from \eqref{ref-3.4-21}
whereas those labeled by $\beta$ are obtained from \eqref{ref-3.5-22}.

The commutativity of this diagram follows easily from lemma
\ref{ref-3.2-18} together with the observation that all
horizontal arrows are obtained by Yoneda multiplying with 
elements of suitable $\Ext$-groups. For example the connecting maps
are obtained from multiplying with the element of $\Ext^1(F_2,F_1)$
representing the exact sequence \eqref{ref-3.6-24}. 

If we now have non-degeneracy for two out of the three pairs
$(F_1,G)$, $(F,G)$, $(F_2,G)$ then we also have it for the third pair
because of the five-lemma.

The case where we have an exact sequence as in
\eqref{ref-3.7-25} is treated similarly. 
\end{proof}
To continue we define a some canonical objects in $\Fscr$. Let
$a\in\NN$. Then we write
\begin{align*}
F^0_{0a}&=(x^{-a}k[x],k,0,\Id_{k[x,x^{-1}]})\\
F^1_{0a}&=(x^{-a}k[x],0,k,\Id_{k[x,x^{-1}]})\\
\end{align*}
(the reason for this notation will become clear  in Section \S\ref{ref-5-36}).
\begin{lemma} 
\label{ref-3.4-26}
Every object $F$ in $\Fscr$ has a finite filtration
$0=F_0\subset \cdots \subset F_n=F$ such that the corresponding subquotients
are among the $F^i_{0a}$.
\end{lemma}
\begin{proof} By the structure of $\Cscr_x$ there must be a surjective
  map $\phi:F_x\r E^{i}$ where $i=0$ or $i=1$. Hence $\im \phi$ is a
  non-trivial quotient. Since it is easy to see that the subobjects of
  $E^i$ in $\Cscr$ are of the form $F^i_{0a}$ we are done.
\end{proof}
Using \eqref{ref-2.4-14} we can compute the $\Hom$ and
$\Ext$-groups between the $F^i_{0a}$. The results are given in the next
lemma. 
\begin{lemma} \label{ref-3.5-27} One has
\begin{align*}
\Hom(F^i_{0a},F^j_{0b})&=
\begin{cases}
k& \text{if $i=j$ and $a\le b$}\\
0&\text{otherwise}
\end{cases}
\\
\Ext^1(F^i_{0a},F^j_{0b})&=
\begin{cases}
k& \text{if $i=1-j$ and $a>b$}\\
0&\text{otherwise}
\end{cases}
\end{align*}
\end{lemma}
\begin{proof}
The claim for $\Hom$ is trivial, so we concentrate on $\Ext$. 

We use \eqref{ref-2.4-14}. This immediately yields that
$\Ext^1(F^i_{0a},F^j_{0b})=0$ if $j\neq 1-i$. If $j=1-i$ then we have the
following exact sequence. 
\begin{equation}
\label{ref-3.8-28}
\Hom_{k[x]}(x^{-a}k[x],x^{-b}k[x])\r k \r
\Ext^1_\Cscr(F^i_{0a},F^j_{0b})\r 0
\end{equation}
This yields that $\Ext^1_\Cscr(F^i_{0a},F^j_{0b})=k$ if and only if
$\Hom_{k[x]}(x^{-a}k[x],x^{-b}k[x])=0$, i.e. if and only if $a>b$.
\end{proof}
We are now in a position to prove the main result of this section.
\begin{theorem} \label{ref-3.6-29} $\Cscr$ satisfies Serre duality.
\end{theorem}
\begin{proof} We show first that $\Cscr$ satisfies Serre duality for objects
$F,G$ in $\Fscr$.
We will show the non-degeneracy of \eqref{ref-3.4-21} and
  \eqref{ref-3.5-22} by induction of $\rk_{k[x]}(F)$,
  $\rk_{k[x]}(G)$. This reduces us to the case where $F=F^i_{0a}$,
  $G=F^j_{0b}$. So we need to check the non-degeneracy of
\begin{align}
  \Hom_{\Cscr}(F^i_{0a},F^j_{0b})\times
  \Ext^1_{\Cscr}(F^j_{0b},F^{i'}_{0,a-1})\r
  \Ext^1_{\Cscr}(F^i_{0a},F^{i'}_{0,a-1})\xrightarrow{}&
  k\label{ref-3.9-30}\\
  \Ext^1_{\Cscr}(F^i_{0a},F^j_{0b})\times
  \Hom_{\Cscr}(F^j_{0b},F^{i'}_{0,a-1})\r
  \Ext^1_{\Cscr}(F^i_{0a},F^{i'}_{0,a-1})\xrightarrow{}& k
 \label{ref-3.10-31}
\end{align}
where $i'=1-i$. We will concentrate ourselves on
\eqref{ref-3.10-31}. \eqref{ref-3.9-30} is similar. By
\eqref{ref-3.5-27} the only non-trivial case is given by
$j=1-i$ and $a>b$. In that case all vector spaces involved are equal to
$k$ and what we want to prove follows from inspecting 
\eqref{ref-3.8-28}.

Now we show that $\Cscr$ has almost split sequences. By Theorem \ref{ref-1.1-0}
this implies that $\Cscr$ satisfies Serre duality.

By Proposition 
\ref{ref-2.2-12} it is clearly sufficient to construct almost split 
sequences ending in indecomposable objects in $\Fscr$ and $\Tscr$.
First let $F\in \Fscr$
be indecomposable. Since $\Ext^1_\Cscr(F,VF)\cong \Hom_\Cscr(F,F)^\ast$, 
$\Ext^1_\Cscr(F,VF)$ has a simple socle as (left or right) $\Hom_\Cscr(F,F)$ module.
Let $\zeta$ be a non-zero element this socle. It is well-known, and easy
to see that $\zeta$ defines the almost split sequence
in $\Fscr$ ending in $F$.
\[
0\r VF \r M\r F \r 0
\]
But this is also an almost split sequence in $\Cscr$ since if $T\in \Tscr$ then
$\Hom(T,F)=0$.

Now let $T\in \Tscr$. Then there exists an almost split sequence in $\Tscr$
\begin{equation}
\label{ref-3.11-32}
0\r T''\r T'\r T\r 0
\end{equation}
(since $\Tscr$ is equivalent to the category of $x$-torsion modules over
$k[x]$). 

We have to prove that any pullback for $C\r T$ of \eqref{ref-3.11-32} with
$C$ indecomposable is split. Clearly we only have to consider the case
$C\in \Fscr$. But then it follows from Proposition \ref{ref-2.2-12} 
that the pullback is split. This finishes the proof.
\end{proof}
\section{Relation with one-dimensional graded type $A_n$-singularities}
\label{ref-4-33}
Let $\Cscr$ be the hereditary category which was described in the
previous section. We will now show that $\Cscr$ can be considered as a
limit of certain graded simple singularities.

If $m\in\NN$ then the graded simple $A_{2m-1}$-singularity of
dimension one is by definition the graded subring $R_m$ of $k[x]\oplus
k[x]$ generated by $u=(x,x)$ and $v=(x^m,0)$. 
It is easy to see that $R\cong k[u,v]/(u^mv-u^{2m})$
and hence this is equivalent to the classical definition (see for example \cite{DW}). We put
$\Cscr_m=\mod(R_m)$. 

Let us consider $k[x]$ as being diagonally embedded in $k[x]\oplus
k[x]$. That is we identify $x$ with $(x,x)$.  Clearly we have
\[
(R_{m})_x=k[x,x^{-1}]\oplus k[x,x^{-1}]
\]
Hence if $M\in \Mod(R_m)$ then $M_x$ is canonically a sum of
two $k[x,x^{-1}]$-modules which we denote by $M^0_x$ and $M^1_x$
respectively. This allows us to define the following functor.
\[
U_m:\Cscr_m\r\Cscr: M\mapsto (M,(M^0_x)_0, (M^1_x)_0,\Id)
\]
Clearly $U_m$ is faithful.
We have inclusions
\[
k[x]\subset \cdots \subset R_{m+1}\subset R_m \subset \cdots
R_0=k[x]\oplus k[x]
\]
Dualizing these yield restriction functors
\[
\Cscr_0\r \cdots \r \Cscr_m\r \Cscr_{m+1}
\r\cdots \r \mod(k[x])
\]
It is clear that these restriction functors are compatible with the
functors $(U_m)_m$. Define $\Cscr_\infty$ as the 2-direct limit of the
$\Cscr_m$. That is the objects in $\Cscr_\infty$ are the objects in
$\coprod_m\Cscr_m$ and we put
\[
\Hom_{\Cscr_\infty}(M,N)=\dirlim \Hom_{\Cscr_m}(M,N)
\]
The functors $(U_m)_m$ define a functor $U_\infty:\Cscr_\infty \r
\Cscr$. 
\begin{proposition}
\label{ref-4.1-34} The functor $U_\infty$ defined above is an
  equivalence. 
\end{proposition}
\begin{proof}
  From the definition it is clear that $U_\infty$ is faithful. So we
  only have to show that it is full and essentially surjective.
  
  We will first show that $U_\infty$ is full. Let $M,N\in \Cscr_m$ and
  let $f:M\r N$ be a homomorphism in $\Cscr$. So $f$ is in fact a
  $k[x]$-linear homomorphism $f:M\r N$ such that the localization 
  $f_x$ is $k[x,x^{-1}]\oplus k[x,x^{-1}]$-linear.

Let $y=(x,0)$. Then $y^n\in R_m$ for $n\ge m$ and $R_n$ as subring of $R_m$ is generated
by $x$ and $y^n$.  To prove fullness of $U_\infty$ it is sufficient
that $f$ is $y^n$-linear for $n\gg 0$. Let $T$ be the torsion submodule of
$N$ and consider the $k[x]$ linear map
$M\r N$ given by $f^{(n)}=f(y^n-)-y^n f(-)$. Since after localizing at $x$, $f$
is $y$ linear, it follows that the image of $f^{(n)}$ lies in $T$.  Since $T$ is rightbounded it is clear
that $f^{(n)}$  must be zero if $n\gg 0$. This proves what we want.

Now we prove essential surjectivity. First let $F\in\Fscr$. Then
we claim that $F\subset F_x$ is stable under multiplication by $y^n$
for $n\gg 0$. First note that $yF$ is a finitely generated $k[x]$-submodule
of $F_x$. Hence $x^n y F\subset F$. Since $x^n y=y^{n+1}$ this proves
what we want. 

Now let $T\in \Tscr$. Then as graded $k[x]$-module $T$ has right
bounded grading and since $k[x]_{<m}=(R_n)_{<m}$ for $n\ge m$ it follows that  for $n\gg 0$
we may consider $T$ as a graded $R_n$-module. This proves what we want. 
\end{proof}
\section{Representation theory}
\label{ref-5-36}

In section we construct the AR-quiver of $\Cscr$. From the above
discussion it follows that the components of the AR-quiver of $\Cscr$
lie either in $\Tscr$ or in $\Fscr$. Since $\Tscr$ is equivalent to the
$x$-torsion modules in $\gr(k[x])$ it has a unique component which is
$\ZZ A_\infty$. So the main difficulty is represented by the component(s)
in $\Fscr$.

We now describe the indecomposable  torsion free
objects in $\Cscr$ as well as the associated Auslander-Reiten quiver
(see \cite{ReVdB1}).  Using Proposition \ref{ref-4.1-34} this
could be easily obtained by using a graded version of the results in
\cite{DW}.  However for completeness we give an independent proof
here.

For $m>0$ denote by $F_{ma}$ the unique indecomposable projective
$R_m$-module in with grading starting in degree $-a$ (thus
$F_{ma}=F_{m0}(a)$).  For $m=0$ we let $F^0_{00}$, $F^1_{00}$ be the
two indecomposable $R_0$ modules whose gradings starts exactly at $0$.
We also put $F^i_{0a}=F^i_{00}(a)$  (as in Section \S\ref{ref-3-15}).

Finally to simplify the notation we will write $F^i_{ma}$
($i=\emptyset$, if $m\neq 0$) for $U_m(F^i_{ma})$.
\begin{proposition} The indecomposable objects in $\Fscr$ are given by 
  $F^i_{ma}$. Furthermore the associated Auslander-Reiten quiver is
  given by Figure \ref{ref-1-35}
\begin{figure}
\begin{center}
$
\xymatrix{
&\ar@{.}[d]&&\ar@{.}[d]&&\ar@{.}[d]&\\       
 \ar@{--}[r] & F_{41}\ar[dr] \ar@{--}[rr] && F_{42}\ar[dr] 
\ar@{--}[rr] && F_{43}\ar@{--}[r]&\\ 
& \ar@{--}[r]& F_{31} \ar[dr]\ar[ur]\ar@{--}[rr] 
&& F_{32}\ar[dr]\ar[ur]  \ar@{--}[r]&&\\ 
 \ar@{--}[r] & F_{20}\ar[ur]\ar[dr] \ar@{--}[rr] && 
F_{21}\ar[ur]\ar[dr] \ar@{--}[rr] 
&& F_{22}\ar@{--}[r]&\\ 
& \ar@{--}[r]& F_{10} \ar[dr]\ar[ddr]\ar[ur]\ar@{--}[rr] 
&& F_{11}\ar[dr]\ar[ddr]\ar[ur]  \ar@{--}[r]&&\\ 
 \ar@{--}[r] & F^0_{0-1}\ar[ur] \ar@{--}[rr] && 
F^1_{00}\ar[ur] \ar@{--}[rr] 
&& F^0_{01}\ar@{--}[r]&\\ 
 \ar@{--}[r] & F^1_{0-1}\ar[uur] \ar@{--}[rr] && 
F^0_{00}\ar[uur] \ar@{--}[rr] 
&& F^1_{01}\ar@{--}[r]&\\ 
}
$
\end{center}
\caption{The Auslander-Reiten quiver of $\Fscr$}
\label{ref-1-35}
\end{figure}
\end{proposition}
\begin{proof}
By Serre duality it follows that $\Ext^1_\Cscr(F^i_{ma},VF^i_{ma})$ is one
dimensional. Therefore its unique (up to scalar multiplication) non-zero
element represents the almost split sequence ending in $F^i_{ma}$. 

Let us now explicitly construct non-split extensions between $F^i_{ma}$ and
$VF^i_{ma}$. First note
\[
VF_{ma}=F_{m,a-1}
\]
where for simplicity we have written $F_{0a}=F^0_{0a}\oplus F^1_{1a}$,
and
\[
V F^i_{0a}=F^{1-i}_{0a-1}
\]
To construct the extension associated to $F_{m,a}$  we note that
$F_{m-1,a-1}$ and $F_{m+1,a}$ are naturally submodules of $F_{m,a}$ whose sum
is $F_{m,a}$ and whose intersection
is $F_{m,a-1}$. Hence the exact sequence
\begin{equation}
\label{ref-5.1-37}
0\r F_{m,a-1}\r F_{m-1,a-1}\oplus F_{m+1,a} \r  F_{m,a}\r 0
\end{equation}
yields the sought extension.

To construct the extension associated to $F^i_{0a}$  we note that
$F_{1a}$ maps surjectively to $F^i_{0a}$ with kernel $F^{1-i}_{0a-1}$. Thus in this
case the sought extension is
\begin{equation}
\label{ref-5.2-38}
0\r F^{1-i}_{0,a-1}\r F_{1a}\oplus F^i_{0a} \r  F_{m,a}\r 0
\end{equation}
It is now easy to assemble the almost split sequences given by \eqref{ref-5.1-37} and
\eqref{ref-5.2-38} into the translation quiver given by Figure 1.

To show that Figure \ref{ref-1-35} is the entire AR-quiver of $\Fscr$ (and not
just a component) we have to show that there are no other
indecomposable objects.

So assume that $F$ is an indecomposable object in $\Fscr$, not
occurring among the $F^i_{ma}$. By lemma \ref{ref-3.4-26} there
exist a non-zero map $F^i_{0a}\r F$ for some $i,a$. Using the defining property of AR-sequences we may use this to construct a non-zero map $F^i_{mb}$ for
some $i$ (possibly $\emptyset$) and $m$, and for $b$ arbitrarily large.  

Now note that the only non-trivial torsion free
quotients of $F^i_{mb}$ are $F^i_{mb}$ itself and $F^{0,1}_{0b}$ (if
$m\neq 0$). Since all these quotients possess a non-trivial element in
degree $-b$ it follows that $\Hom_\Cscr(F^i_{mb},F)=0$ for $b\gg 0$. This
finishes the proof.
\end{proof}
\def\cprime{$'$} \def\cprime{$'$}
\ifx\undefined\bysame
\newcommand{\bysame}{\leavevmode\hbox to3em{\hrulefill}\,}
\fi

\end{document}